\documentclass[amstex,12pt, amssymb]{article}

\usepackage{mathtext}
\usepackage[cp1251]{inputenc}
\usepackage[T2A]{fontenc}
\usepackage[dvips]{graphicx}
\usepackage{amsmath}
\usepackage{amssymb}
\usepackage{amsxtra}
\usepackage{latexsym}
\usepackage{ifthen}

\newcounter{lemma}[section]

\newcounter{corollary}[section]

\newcounter{remark}[section]

\newcounter{theorem}[section]

\newcounter{proposition}[section]

\numberwithin{equation}{section}

\def\XXint#1#2#3{{\setbox0=\hbox{$#1{#2#3}{\int}$}
     \vcenter{\hbox{$#2#3$}}\kern-.5\wd0}}

\def\cc{\setcounter{equation}{0}
\setcounter{figure}{0}\setcounter{table}{0}}

\begin{document}

\markboth{\centerline{VLADIMIR RYAZANOV}} {\centerline{STIELTJES
INTEGRALS FOR HARMONIC AND ANALYTIC FUNCTIONS}}

\author{{Vladimir Ryazanov}}

\title{{\bf The Stieltjes integrals in the theory\\ of harmonic functions}}

\maketitle

\large \begin{abstract} We study various Stieltjes integrals as
Poisson-Stieltjes, conjugate Poisson-Stieltjes, Schwartz-Stieltjes
and Cauchy-Stieltjes and prove theorems on the existence of their
finite angular limits a.e. in terms of the Hilbert-Stieltjes
integral. These results hold for arbitrary bounded integrands that
are differentiable a.e. and, in particular, for integrands of the
class ${\cal {CBV}}$ (countably bounded variation).
\end{abstract}


\large
\cc
\section{Introduction}

\medskip

Recall that a path in $\mathbb D:=\{ z\in\Bbb C : |z|<1\}$
terminating at $\zeta\in\partial\mathbb D$ is called {\bf
nontangential} if its part in a neighborhood of $\zeta$ lies inside
of an angle in $\mathbb D$ with the vertex at $\zeta$. Hence the
limit along all nontangential paths at $\zeta\in\partial\mathbb D$
also named {\bf angular} at the point. The latter is a traditional
tool of the geometric function theory, see e.g. monographs
\cite{Du}, \cite{Ko}, \cite{L}, \cite{Po} and \cite{P}.

It was proved in the previous paper \cite{R3}, see also \cite{R4},
that a harmonic function $u$ given in the unit disk $\mathbb D$ of
the complex plane $\mathbb C$ has angular limits at a.e. point
$\zeta\in\partial\mathbb D$ if and only if its conjugate harmonic
function $v$ in $\mathbb D$ is so. This is the key fact together
with Lemma 2 in Section 5 further to establish the existence of the
so--called Hilbert-Stieltjes integral for a.e.
$\zeta\in\partial\mathbb D$ and the corresponding result on the
angular limits of Cauchy-Stieltjes integral under fairly general
assumptions on integrands, cf. e.g. \cite{H}, \cite{Pr} and
\cite{Sm}, see also \cite{Dy}, \cite{Ki}, \cite{Ko} and \cite{P}.


It is known the very delicate fact due to Lusin that harmonic
functions in the unit circle with continuous (even absolutely
continuous !) boundary data can have conjugate harmonic functions
whose boundary data are not continuous functions, furthemore, they
can be even not essentially bounded in neighborhoods of each point
of the unit circle, see e.g. Theorem VIII.13.1 in \cite{Bari}. Thus,
a correlation between boundary data of conjugate harmonic functions
is not a simple matter, see also I.E in \cite{Ko}.


Denote by $h^p$, $p\in(0,\infty)$, the class of all harmonic
functions $u$ in $\mathbb D$ with bounded $L^p-$norms over the
circles $|z|=r\in(0,1)$. It is clear that $h^p\subseteq
h^{p^{\prime}}$ for all $p>p^{\prime}$ and, in particular,
$h^p\subseteq h^{1}$ for all $p>1$. It is significant that every
function in the class $h^1$ has a.e. nontangential boundary limits,
see e.g. Corollary IX.2.2 in \cite{Go}. It is also known that a
harmonic function $u$ in $\Bbb D$ can be represented as the Poisson
integral
\begin{equation}\label{Poisson}
u(re^{i\vartheta})\ =\ \frac{1}{2\pi}\
\int\limits_{-\pi}\limits^{\pi}\frac{1-r^2}{1-2r\cos(\vartheta-t)+r^2}\
\varphi(t)\ dt
\end{equation} with a function $\varphi\in L^p(-\pi,\pi),\: p>1$, if and only if $u\in
h^p$, see e.g. Theorem IX.2.3 in \cite{Go}. Thus, $u(z)\to
\varphi(\vartheta)$ as $z\to e^{i\vartheta}$ along any nontangential
path for a.e. $\vartheta$, see e.g. Corollary IX.1.1 in \cite{Go}.
Moreover, $u(z)\to \varphi(\vartheta_0)$ as $z\to e^{i\vartheta_0}$
at points $\vartheta_0$ of continuity of the function $\varphi$, see
e.g. Theorem IX.1.1 in \cite{Go}.


Note also that $v\in h^p$ whenever $u\in h^p$ for all $p>1$ by the
M. Riesz theorem, see \cite{R}, see also Theorem IX.2.4 in
\cite{Go}. Generally speaking, this fact is not trivial but it
follows immediately for $p=2$ from the Parseval equality, see e.g.
the proof of Theorem IX.2.4 in \cite{Go}. The case $u\in h^1$ is
more complicated.


As known, a harmonic function $u$ in $\Bbb D$ belongs to $h^1$ if
and only if it can be represented as the Poisson--Stieltjes integral
\begin{equation}\label{Poisson-Stieltjes}
u(re^{i\vartheta})\ =\ \frac{1}{2\pi}\
\int\limits_{-\pi}\limits^{\pi}\frac{1-r^2}{1-2r\cos(\vartheta-t)+r^2}\
d \Phi(t)
\end{equation} with $\Phi : [-\pi,\pi]\to\mathbb R$ of bounded variation, see e.g.
Theorem IX.2.2 in \cite{Go}. Moreover, if the function $\Phi$ has a
finite derivative at a point $\vartheta_0\in(-\pi,\pi)$, then
$u(z)\to\Phi^{\prime}(\vartheta_0)$ as
$z:=re^{i\vartheta}\to\zeta_0:=e^{i\vartheta_0}$ along all
nontangential paths in $\mathbb D$ to the point $\zeta_0$, see e.g.
Theorem IX.1.4 in \cite{Go}.

\medskip

The present paper is devoted to the study of the corresponding
Stieltjes integrals in the case when $\Phi$ is generally speaking
not of bounded variation. The article has a methodological
character. Namely, in Section 2 the first base is that the formula
of integration by parts holds with no regularity conditions on
functions if at least one of the integrals exists, see Lemma 1, and
that the latter is true if one of the functions is absolute
continuous and the second one is bounded and its set of points of
discontinuity has measure zero, see Proposition 1. Moreover, it is
demonstrated by Example 1 that the boundedness condition is
essential. Finally, by Remark 1 measure zero for the set of points
of discontinuity is also necessary. Just thanking to this
methodology basis, it has been possible to extend all the rest
results on various Stieltjes integrals to the case of absence of
bounded variation.

\medskip

\cc

\section{Expansion of the Riemann--Stieltjes integral}
First of all, recall a classical definition of the
Riemann--Stieltjes integral. Namely, let $\rm I = [a, b]$ be a
compact interval in $\Bbb R.$ A {\bf partition} $\rm P$ of $\rm I$
is a collection of points $t_0,t_1,\ldots ,t_p\in \rm I$ such that
$a = t_0 \le t_1 \le \ldots \le t_p = b$. Now, let $g: \rm I\to\Bbb
R$ and $f: \rm I\to\Bbb R$ be bounded functions and let $f$ be in
addition nondecreasing. The {\bf Riemann--Stieltjes integral} of $g$
with respect to $f$ is a real number $A$, written
$A=\int\limits_{\rm I}g\, df$, if, for every $\varepsilon >0$, there
is $\delta >0$ such that
$$
|\ \sum\limits_{k=1}^p g(\tau_k)[f(t_k)-f(t_{k-1})]\ -\ A\ |\ <\
\varepsilon
$$
for every partition $\rm P=\{ t_0,t_1,\ldots ,t_p\} $ of $\rm I$
with $|t_k-t_{k-1}|\le\delta$, $k=1,\ldots , p$, and every
collection $\tau_k\in[t_{k-1},t_k]$, $k=1,\ldots , p$. In other
words and notations, \begin{equation}\label{eqLIM} \int\limits_a^b\,
g\, df\ :=\ \lim\limits_{\delta\to 0}\ \sum\limits_{k=1}^p\
g(\tau_k)\cdot \Delta_k f\ \ \ \ \ \ \ \hbox{as}\ \ \ \delta :=
\max\limits_{k=1,\ldots , p}\, |t_k-t_{k-1}|\to 0
\end{equation}
where $\Delta_k f:=f(t_k)-f(t_{k-1})$, $k=1,\ldots , p$, if a finite
limit in (\ref{eqLIM}) exists and it is uniform with respect to
partitions $\{ t_k\}$ and intermediate points $\{\tau_k\}$.

\medskip

We extend the definition of the Riemann--Stieltjes integral to
arbitrary functions $g$ and $f$ for which the limit (\ref{eqLIM})
exists. Let us start from the following general fact, cf. e.g.
\cite{Bo}, \cite{HT}, \cite{Pf} and \cite{S}.

\medskip

{\bf Lemma 1.} {\it Let $\rm I=[a,b]$, $g: \rm I\to\Bbb R$ and $f:
\rm I\to\Bbb R$ be arbitrary functions. If one of the integrals
$\int g\, df$ and $\int f\, dg$ exists, then the second one is so
and
\begin{equation}\label{eqPARTS} \int\limits_{\rm I}\,
g\, df\ +\ \int\limits_{\rm I}\, f\, dg\ =\ g(b)\cdot f(b)\ -\
g(a)\cdot f(a)\ .
\end{equation}}

\begin{proof}
For definiteness, let us assume that there exists the integral $\int
f\, dg$. Then the identity for the arbitrary integral sum of the
integral $\int g\, df$:
$$
\sum\limits_{k=1}^pg(\tau_k)\cdot(f(t_k)-f(t_{k-1})) = -g(t_0)f(t_0)
- \sum\limits_{k=1}^{p+1} f(t_k)\cdot(g(\tau_k)-g(\tau_{k-1})) +
g(t_p)f(t_p),
$$
where $\tau_0:=t_0$ and $\tau_{p+1}=t_p$, implies the desired
conclusions.
\end{proof} $\Box$

\bigskip

By Theorem 13.1.b in \cite{HT}, we have the significant consequence
of Lemma 1.

\medskip

{\bf Proposition 1.}\label{pr1} {\it Let $\rm I=[a,b]$, a function
$f: \rm I\to\Bbb R$ be absolutely continuous and let $g: \rm
I\to\Bbb R$ be a bounded function whose set of points of
discontinuity is of measure zero. Then both integrals
$\int\limits_{\rm I} g\, df$ and $\int\limits_{\rm I} f\, dg$ exist
and the relation (\ref{eqPARTS}) holds.}

\medskip

{\bf Example 1.} {\it The condition that the function $g$ is bounded
in Proposition 1 is essential. It is clear from the simplest example
of the pair of the functions on $[0,1]$: $f(t)=t$ and $g(t)=0$
except the points $t_n=1/n$ where $g(t_n):=n^2$, $n=1,2,\ldots $.
Indeed, we see that the lower limit of the integral sums is $0$ and
the upper limit is $\infty$.}

\medskip

{\bf Remark 1.} The condition on measure zero for the set of points
of discontinuity of $g$ in Proposition 1 is also necessary. Indeed,
take the case $f(t)\equiv t$ and consider a function $g$ with its
set $S$ of points of discontinuity of a length $l>0$. By
subadditivity $l\le \sum l_n$ where $l_n$, $n=1,2,\ldots $, is the
length of the set $S_n$ of points in $S$ with jumps which are
greater or equal to $1/n$. Hence $l_N>0$ for some $N=1,2,\ldots $.
Let us cover every point $t\in S_N\setminus\{ a,b\}$ by all
intervals in $(a,b)$ centered at $t$ whose lengths are less than an
arbitrary prescribed $\delta>0$. By the Vitali theorem there is a
countable collection of such mutually disjoint intervals covering
almost every point of $S_N$, see e.g. Theorem IV.3.1 in \cite{S}.
The sum of their lengths is not less than $l_N$ and all such
intervals contain jumps of $g$ that are not less than $1/N$. Thus,
the difference between lower and upper limits of integrals sums is
not less than $N^{-1}l_N$, i.e., it is not zero.

It is clear that formula (\ref{eqPARTS}) is also valid for complex
va\-lu\-ed functions of the natural parameter on rectifiable Jordan
curves because their real and imaginary parts can be considered as
real valued functions on segments of $\mathbb R$. Moreover, the
corresponding statements hold on closed Jordan curves $\rm J$ with
the relation
\begin{equation}\label{eqJORDAN} \int\limits_{\rm J}\,
g\, df\ =\ -\int\limits_{\rm J}\, f\, dg\
\end{equation}
where we should apply {\bf cyclic partitions} $\rm P$ of $\rm J$ by
collections of cyclic ordered points $\zeta_0,\zeta_1, \ldots ,
\zeta_p$ on $\rm J$ with $\zeta_0=\zeta_p$.

\cc
\section{On the Poisson--Stieltjes integrals}

Recall that the {\bf Poisson kernel} is the $2\pi -$periodic
function
\begin{equation}\label{eqP} P_r(\Theta)\ =\
\frac{1-r^2}{1-2r\cos\Theta +r^2}\ , \ r<1\ , \ \Theta\in\mathbb R\
.
\end{equation}
By Proposition 1 and Remark 1, the {\bf Poisson-Stieltjes integral}
\begin{equation}\label{eqPS} \mathbb U(z)\ =\ \mathbb U_{\Phi}(z)\ :=\
\frac{1}{2\pi}\ \int\limits_{-\pi}\limits^{\pi}P_r(\vartheta -t)\
d\Phi(t)\ ,\ \ \ z=re^{i\vartheta}, \ r<1\ ,\ \vartheta\in\mathbb R\
\end{equation}
is well-defined for $2\pi -$periodic continuous functions,
furthermore, for bounded functions  $\Phi:\mathbb R\to\mathbb R$
whose set of points of discontinuity is of measure zero because the
function $P_r(\Theta)$ is continuously differentiable and hence it
is absolutely continuous.

Moreover, directly by the definition of the Riemann-Stieltjes
integral and the Weierstrass type theorem for harmonic functions,
see e.g. Theorem I.3.1 in \cite{Go}, $U$ is a harmonic function in
the unit disk $\mathbb D:=\{ z\in\mathbb C: |z|<1\}$ because the
function $P_r(\vartheta -t)$ is the real part of the analytic
function
\begin{equation}\label{eqA}
{\cal A}_{\zeta}(z)\ :=\ \frac{\zeta +z}{\zeta -z}\ ,\ \ \
\zeta=e^{it},\ \ \ z=re^{i\vartheta}\ , \ r<1\ ,\ \vartheta \ {\rm
and}\ t\in\mathbb R\ . \end{equation}

\bigskip

{\bf Theorem 1.} {\it Let $\Phi:\mathbb R\to\mathbb R$ be a  $2\pi
-$periodic bounded function whose set of points of discontinuity has
measure zero. Suppose that $\Phi$ is differentiable at a point
$t_0\in\mathbb R$. Then
\begin{equation}\label{eqLP}
\lim\limits_{z\to\zeta_0}\ \mathbb U_{\Phi}(z)\ =\
\Phi^{\prime}(t_0)
\end{equation} along all
nontangential paths in $\mathbb D$ to the point
$\zeta_0:=e^{it_0}\in\partial\mathbb D$.}

\bigskip

\begin{proof} Indeed, by Proposition 1 with $g(t):=\Phi(t)$ and $f(t):=P_r(\vartheta
-t)$, $t\in\mathbb R$, for every fixed $z=re^{i\vartheta}$, $r<1$,
$\vartheta\in\mathbb R$, we obtain that
\begin{equation}\label{eqPS2}
\int\limits_{-\pi}\limits^{\pi} P_r(\vartheta -t)\ d\Phi(t)\ =\
\int\limits_{-\pi}\limits^{\pi} \Phi(t)\cdot
\frac{\partial}{\partial\vartheta}\, P_r(\vartheta -t)\ dt\ \ \ \
\forall\  r\in(0,1)\ ,\ \vartheta\in\mathbb R\ ,
 \end{equation}
because of the $2\pi -$periodicity of the given functions $g$ and
$f$ the right hand side in (\ref{eqPARTS}) is equal to zero, $f\in
C^1$ and
$$
\frac{\partial}{\partial\vartheta}\, P_r(\vartheta -t)\ =\ -\
\frac{\partial}{\partial t}\, P_r(\vartheta -t)\ \ \ \ \forall\
r\in(0,1)\ ,\ \vartheta \ {\rm and}\ t\in\mathbb R\ .
$$

Now, considering the Poisson integral
$$
u(re^{i\vartheta})\ :=\ \frac{1}{2\pi}\
\int\limits_{-\pi}\limits^{\pi}P_r(\vartheta - t)\ \Phi(t)\ dt\ ,
$$
we see by the Fatou result, see e.g. 3.441 in \cite{Z}, p. 53, or
Theorem IX.1.2 in \cite{Go}, that
$\frac{\partial}{\partial\vartheta}\ u(z)\to \Phi^{\prime}(t_0)$ as
$z\to\zeta_0$ along any nontangential path in $\mathbb D$ ending at
$\zeta_0$. Thus, the conclusion follows because just $\mathbb
U_{\Phi}(z)=\frac{\partial}{\partial\vartheta}\ u(z)$  by
(\ref{eqPS2}).
\end{proof} $\Box$

\bigskip

{\bf Corollary 1.} {\it   If $\Phi:\mathbb R\to\mathbb R$ is a $2\pi
-$periodic continuous differentiable a.e. function, then $\mathbb
U_{\Phi}(z)\to\Phi^{\prime}(\arg\zeta)$ as $z\to\zeta$ for a.e.
$\zeta\in\partial\mathbb D$ along all nontangential paths in
$\mathbb D$ to the point $\zeta$.}

\bigskip

Here we denote by $\arg \zeta$ the {\bf principal branch of the
argument} of $\zeta\in\mathbb C$ with $|\zeta |=1$, i.e., the unique
number $\tau\in (-\pi ,\pi ]$ such that $\zeta =e^{i\tau}.$

\bigskip

{\bf Remark 2.} Note that the function of interval
$\Phi_*([a,b]):=\Phi(b)-\Phi(a)$ ge\-ne\-ral\-ly speaking generates
no finite signed measure (charge) if $\Phi$ is not of bounded
variation. Hence we cannot apply the known Fatou result on the
angular boun\-da\-ry limits directly to the Poisson--Stieltjes
integrals, see e.g. Theorem I.D.3 in \cite{Ko}.

\bigskip

{\bf Corollary 2.} {\it   If $\Phi:\mathbb R\to\mathbb R$ is a $2\pi
-$periodic bounded function that is differentiable a.e., then
$\mathbb U_{\Phi}(z)\to\Phi^{\prime}(\arg\zeta)$ as $z\to\zeta$ for
a.e. $\zeta\in\partial\mathbb D$ along all nontangential paths in
$\mathbb D$ to the point $\zeta$.}

\bigskip

We call $\Phi:\mathbb R\to\mathbb C$  a function of {\bf countably
bounded variation} and write $\Phi\in\mathcal{CBV}(\mathbb R)$ if
there is a countable collection of mutually disjoint intervals
$(a_n,b_n)$, $n=1,2,\ldots$ on each of which the restriction of
$\Phi$ is of bounded variation and the set $\mathbb R\setminus
\bigcup\limits_1\limits^{\infty}(a_n,b_n)$ is countable.

\bigskip

{\bf Corollary 3.} {\it   If $\Phi:\mathbb R\to\mathbb R$ is a $2\pi
-$periodic bounded function of the class $\mathcal{CBV}(\mathbb R)$,
then $\mathbb U_{\Phi}(z)\to\Phi^{\prime}(\arg\zeta)$ as $z\to\zeta$
for a.e. $\zeta\in\partial\mathbb D$ along all nontangential paths
in $\mathbb D$ to the point $\zeta$.}

\cc
\section{On the conjugate Poisson--Stieltjes integrals}

Recall that the {\bf conjugate Poisson kernel} is the $2\pi
-$periodic function
\begin{equation}\label{eqQ} Q_r(\Theta)\ =\
\frac{2r\sin\Theta}{1-2r\cos\Theta +r^2}\ , \ r<1\ , \
\Theta\in\mathbb R\ .
\end{equation}
By Proposition 1 the {\bf conjugate Poisson-Stieltjes integral}
\begin{equation}\label{eqCPS} \mathbb V(z)\ =\ \mathbb V_{\Phi}(z)\ :=\
\frac{1}{2\pi}\ \int\limits_{-\pi}\limits^{\pi}Q_r(\vartheta -t)\
d\Phi(t)\ ,\ \ \ z=re^{i\vartheta}, \ r<1\ ,\ \vartheta\in\mathbb R\
,
 \end{equation}
is well-defined for $2\pi -$periodic bounded functions $\Phi:\mathbb
R\to\mathbb R$ whose set of points of discontinuity is of measure
zero because the function $Q_r(\Theta)$ is con\-ti\-nuous\-ly
differentiable and hence it is absolutely continuous. Again,
directly by the definition of the Riemann-Stieltjes integral and the
Weierstrass type theorem $\mathbb V_{\Phi}$ is a conjugate harmonic
function for $\mathbb U_{\Phi}$ in the unit disk $\mathbb D$ because
the function $Q_r(\vartheta -t)$ is the imaginary part of the same
analytic function (\ref{eqA}).

\bigskip

By Theorem 1 in \cite{R3}, see also \cite{R4}, we have the following
significant consequences from Theorem 1 and Corollaries 1--3.

\bigskip

{\bf Corollary 4.} {\it Let $\Phi:\mathbb R\to\mathbb R$ be a  $2\pi
-$periodic continuous function that is differentiable a.e. Then
$\mathbb V_{\Phi}(z)$ has a finite limit $\varphi(\zeta)$ as
$z\to\zeta$ along all nontangential paths in $\mathbb D$ to a.e.
$\zeta\in\partial\mathbb D$.}

\bigskip

{\bf Corollary 5.} {\it Let $\Phi:\mathbb R\to\mathbb R$ be a  $2\pi
-$periodic bounded function that is differentiable a.e. Then
$\mathbb V_{\Phi}(z)$ has a finite limit $\varphi(\zeta)$ as
$z\to\zeta$ along all nontangential paths in $\mathbb D$ to a.e.
$\zeta\in\partial\mathbb D$.}

\bigskip

{\bf Corollary 6.} {\it Let $\Phi:\mathbb R\to\mathbb R$ be a $2\pi
-$periodic bounded function of the class $\mathcal{CBV}(\mathbb R)$.
Then $\mathbb V_{\Phi}(z)$ has a finite angular limit
$\varphi(\zeta)$ as $z\to\zeta$ for a.e. $\zeta\in\partial\mathbb
D$.}

\bigskip

The function $\varphi(\zeta)$ will be calculated in the explicit
form through $\Phi(\zeta)$ in terms of the so--called
Hilbert--Stieltjes integral. To prove this fact we need first to
establish one auxiliary result in the next section.

\bigskip

\cc
\section{On the Hilbert--Stieltjes integral}

{\bf Lemma 2.} {\it Let $\Phi:\mathbb R\to\mathbb R$ be a $2\pi
-$periodic bounded function whose set of points of discontinuity has
measure zero. Suppose that $\Phi$ is differentiable at a point
$t_0\in\mathbb R$. Then the difference
\begin{equation}\label{eqPL}
\mathbb V_{\Phi}(z)\ -\ \frac{1}{2\pi}\int\limits_{1-|z| \leqslant
|t_0-t| \leqslant \pi} \frac{d\, \Phi(t)}{\tan \frac{t_0-t}{2}}
\end{equation}
converges to zero as $z\to\zeta_0:=e^{it_0}\in\partial\mathbb D$
along the radius in $\mathbb D$ to the point $\zeta_0$.}

\bigskip

\begin{proof} First of all, in the case of need  applying simultaneous rotations
$\zeta_0$ to $\zeta =e^{it}\in\partial\mathbb D$ and $z\in\mathbb D$
in (\ref{eqA}), we may assume that $t_0=0$. Moreover, with no loss
of generality we may assume that $\Phi(0)=0$ and
$\Phi^{\prime}(0)=0$ because, for the linear function
$\Phi_*(t):=\Phi(0)+\Phi^{\prime}(0)\cdot t : (-\pi,\pi]\to\mathbb
R$ extended $2\pi -$periodically to $\mathbb R$, $d\,\Phi_*(t)\equiv
\Phi^{\prime}(0)\, dt$, gives identical zero in the difference
(\ref{eqPL}) in view of the oddness of the kernel $Q_r$ and
$\tan\frac{t}{2}$.

Note that by the oddness of $Q_r$ we have also that $$\mathbb
V_{\Phi}(r)\ =\ \frac{1}{2\pi}\
\int\limits_{-\pi}\limits^{\pi}Q_r(-t)\ d\,\Phi(t)\ =\
\frac{1}{2\pi}\ \int\limits_{-\pi}\limits^{\pi}Q_r(t)\ d\,\Phi(-t)\
,\ \ \forall\ r\in(0,1)\ .
$$
Then the difference (\ref{eqPL}) is split into two parts with
$\varepsilon =\varepsilon(r):=1-r$ :
$$
{\rm I}\ :=\ \frac{1}{2\pi}\
\int\limits_{-\varepsilon}\limits^{\varepsilon}Q_r(t)\ d\,\Phi(-t)\
, \ \ \ \ \ {\rm II}\ :=\ \frac{1}{2\pi}\
\int\limits_{\varepsilon\le t\le\pi}\{ Q_r(t)-Q_1(t)\}\ d\,\{
\Phi(-t)-\Phi(t)\}\ .
$$

Integrating I by parts, we have by  Proposition 1 and Remark 1 and
the oddness of $Q_r(t)$
\begin{equation}\label{eqIP}
{\rm I}\ =\ \frac{1}{2\pi}\
Q_r(\varepsilon)\{\Phi(-\varepsilon)-\Phi(\varepsilon)\}\ +\
\frac{1}{\pi}\
\int\limits_{0}\limits^{\varepsilon}\{\Phi(-t)-\Phi(t)\}\ d\,
Q_r(t)\ .
\end{equation} The first summand converges to zero as
$\varepsilon\to 0$ because $|\Phi(\pm\varepsilon)|=o(\varepsilon)$
and
\begin{equation}\label{eqE}
Q_r(\varepsilon)\ =\ \frac{2r\sin \varepsilon}{1-2r\cos \varepsilon
+r^2}\ =\ \frac{2r\sin \varepsilon}{\varepsilon^2 +
4r\sin^2\frac{\varepsilon}{2}} \le\ 2\,\frac{\sin
\varepsilon}{\varepsilon^2}\ \le\ \frac{2}{\varepsilon}\ .
\end{equation} To estimate the second summand in (\ref{eqIP}) note that $\sin^2
\frac{t}{2}\le \left[\frac{1-r}{2}\right]^2$ and, thus,
$$
Q^{\prime}_r(t)\ =\ \frac{2r\cos t}{1-2r\cos t +r^2}\ -\
\frac{4r^2\sin^2 t}{(1-2r\cos t +r^2)^2}\ =
$$
$$
=\ 2r\ \frac{(1+r^2)\cos t - 2r}{(1-2r\cos t +r^2)^2}\ =\ 2r\
\frac{(1-r)^2-2(1+r^2)\sin^2 \frac{t}{2}}{(1-2r\cos t +r^2)^2}\ \ge
$$
$$
\ge\ 2r\ \frac{(1-r)^2[1-(1+r^2)/2]}{(1-2r\cos t +r^2)^2}\ =\
\frac{r(1+r)(1-r)^3}{(1-2r\cos t +r^2)^2}\ ,
$$
i.e., $Q_r^{\prime}(t)>0$ for all $t\in [ 0 ,\varepsilon]$. Since
$Q_r(t)$ is smooth, it is strictly increasing on $[ 0,\varepsilon]$.
Hence the modulus of the second summand has the upper bound
$$
\frac{1}{\pi}\cdot Q_r(\varepsilon)\cdot \sup\limits_{t\in[ 0,
\varepsilon]}\ \{|\Phi(-t)|+|\Phi(t)|\}\ \le\
\frac{1}{\pi}\cdot\frac{2}{\varepsilon}\cdot o(\varepsilon)\ =\ o(1)
$$
where the inequality follows by (\ref{eqE}). Thus, the second
summand in (\ref{eqIP}) also converges to zero as $\varepsilon\to
0$.

Now, by oddness of the kernels $Q_r(t)$, $r\in(0,1)$ and $Q_1(t)$ we
obtain that
$$
{\rm II}\ :=\ \frac{1}{\pi}\
\int\limits_{\varepsilon}\limits^{\pi}\{ Q_r(t)-Q_1(t)\}\ d\,\{
\Phi(-t)-\Phi(t)\}
$$
where
$$
Q_1(t)-Q_r(t)\ =\ \frac{2\sin t}{2(1-\cos t)}\ -\ \frac{2r\sin
t}{\varepsilon^2+2r(1-\cos t)}\ =
$$
$$=\ \frac{2\sin
t}{4\sin^2\frac{t}{2}}\ -\ \frac{2r\sin
t}{\varepsilon^2+4r\sin^2\frac{t}{2}} \ =\ \frac{2\varepsilon^2\sin
t}{4\left(\varepsilon^2+4r\sin^2\frac{t}{2}\right)\sin^2\frac{t}{2}}
\ .
$$
Integrating by parts, we see that the latter integral is equal to
$$
\{ Q_1(\varepsilon)-Q_r(\varepsilon)\}\cdot \{
\Phi(-\varepsilon)-\Phi(\varepsilon)\}\ +\
\int\limits_{\varepsilon}\limits^{\pi}\ \{ \Phi(-t)-\Phi(t)\}\ d\,
\{ Q_1(t)-Q_r(t)\}
$$
Here the first summand converges to zero because
$\Phi(-\varepsilon)-\Phi(\varepsilon)=o(\varepsilon)$ and
\begin{equation}  \label{eqEE}
Q_1(\varepsilon)-Q_r(\varepsilon)\ =\ \frac{2\varepsilon^2\sin
\varepsilon}{4(\varepsilon^2+4r\sin^2\frac{\varepsilon}{2})\sin^2\frac{\varepsilon}{2}}\
\sim\ \frac{1}{\varepsilon}\ \ \ \mbox{as}\ \ \ \varepsilon\to 0\ .
\end{equation}

Thus, it remains to estimate the integral
$$
{\rm III}\ :=\ \int\limits_{\varepsilon}\limits^{\pi}\ \{
\Phi(-t)-\Phi(t)\}\ d\, \{ Q_1(t)-Q_r(t)\}\ =\
\int\limits_{\varepsilon}\limits^{\pi}\ \varphi(t)\ d\, \alpha_r(t)
$$
where $\varphi(t)=\Phi(-t)-\Phi(t)$ and $\alpha_r(t)=Q_1(t)-Q_r(t)$.
To make it first of all note that $\alpha^{\prime}_r(t)<0$ for
$t\in(\varepsilon,\pi)$ because of
$$
\alpha^{\prime}_r(t)\ =\ \frac{2\varepsilon^2\cos
t}{4(\varepsilon^2+4r\sin^2\frac{t}{2})\sin^2\frac{t}{2}}\ -\
                            \frac{2\varepsilon^2\sin^2t\
                            (\varepsilon^2+8r\sin^2\frac{t}{2})}{4[(\varepsilon^2+4r\sin^2\frac{t}{2})\sin^2\frac{t}{2}]^2}\
                            =\
$$
$$
=\ 2\cdot\frac{\varepsilon^2}{\delta^2}\cdot \left[\cos
t\cdot\sin^2\frac{t}{2}\cdot\left(\varepsilon^2\ +\
4r\sin^2\frac{t}{2}\right)\ -\ \sin^2t\cdot\left(\varepsilon^2\ +\
8r\sin^2\frac{t}{2}\right)\right]\ =\
$$
$$
=2\cdot\frac{\varepsilon^2}{\delta^2} \cdot \left[
\varepsilon^2\cdot   \left( \cos t\cdot\sin^2\frac{t}{2} -
\sin^2t\right) - 4r\sin^2\frac{t}{2}\cdot\left( 2\sin^2t - \cos
t\cdot\sin^2\frac{t}{2}\right) \right] =
$$
$$
=2\cdot\frac{\varepsilon^2}{\delta^2} \cdot \left[
\varepsilon^2\cdot\sin^2\frac{t}{2}\cdot \left( \cos t -
4\cos^2\frac{t}{2}\right) - 4r\sin^4\frac{t}{2}\cdot\left(
8\cos^2\frac{t}{2} - \cos t\right) \right] =
$$
$$
= -2\cdot\frac{\varepsilon^2}{\delta^2}\cdot\sin^2\frac{t}{2}\cdot
\left[ \varepsilon^2\cdot \left( 1+ 2\cos^2 \frac{t}{2} \right) +
4r\sin^2\frac{t}{2}\cdot\left( 1+6\cos^2 \frac{t}{2} \right) \right]
$$
where we many times applied the trigonometric identities $\sin
t=2\sin\frac{t}{2}\cos\frac{t}{2}$ and $1-\cos
t=2\sin^2\frac{t}{2}$, $1+\cos t=2\cos^2\frac{t}{2}$, and the
notation
$$
\delta\ :=\
{2\sin^2\frac{t}{2}\left(\varepsilon^2+4r\sin^2\frac{t}{2}\right)}\
.
$$
The above expression for $\alpha^{\prime}_r(t)$ also implies that
$|\alpha^{\prime}_r(t)|\le c\cdot\frac{\varepsilon^2}{t^4}$. Thus,
$$
|{\rm III}|\ \le\ c\cdot\varepsilon^2
\int\limits_{\varepsilon}\limits^{\pi}|\varphi(t)|\ \frac{d\,
t}{t^4}\ .
$$
Let us fix an arbitrary $\epsilon>0$ and choose a small enough
$\eta>0$ such that $|\varphi(t)|/t<\epsilon$ for all $t\in(0,\eta)$.
Note that we may assume here that $\varepsilon<\eta^2\surd\epsilon$
for small enough $\varepsilon$. Consequently, we have the following
estimates
$$
|{\rm III}|\ \le\ c\cdot\varepsilon^2\cdot\epsilon
\int\limits_{\varepsilon}\limits^{\eta} \frac{d\, t}{t^3}\ +\
c\cdot\varepsilon^2 \int\limits_{\eta}\limits^{\pi}|\varphi(t)|\
\frac{d\, t}{t^4}\ \le\ \frac{c}{2}\cdot\epsilon\ +\
c\pi\cdot\epsilon\cdot M
$$
where $M=\sup\limits_{t\in[0,\pi]}|\varphi(t)|$. In view of
arbitrariness of $\varepsilon$ and $\epsilon$, we conclude that the
integral ${\rm III}$ converges to zero as $\varepsilon \to 0$.
\end{proof} $\Box$

\bigskip

{\bf Theorem 2.} {\it Let $\Phi:\mathbb R\to\mathbb R$ be a $2\pi
-$periodic bounded function. Suppose that $\Phi$ is differentiable
a.e. Then
\begin{equation}   \label{eqLPH}
\lim\limits_{z\to\xi}\ \mathbb V_{\Phi}(z)\ =\
\frac{1}{2\pi}\int\limits_{-\pi}^{\pi} \frac{d\, \Phi(t)}{\tan
\frac{\tau-t}{2}}\ ,\ \ \ \ \ \ \xi:=e^{i\tau}\in\partial\mathbb
 D\ ,
\end{equation}
for a.e. $\tau\in\mathbb R$ along all nontangential paths in
$\mathbb D$ to the point $\xi$.}

\bigskip

Here the singular integral from the right hand side in (\ref{eqLPH})
is understood as a limit of the corresponding proper integrals ({\bf
principal value by Cauchy}):
\begin{equation}\label{eqI}
\mathbb H_{\Phi}(\tau) := \frac{1}{2\pi}\lim\limits_{\varepsilon\to
+0}\ \int\limits_{\varepsilon\le|\tau -t|\le \pi}
\frac{d\,\Phi(t)}{\tan \frac{\tau -t}{2}}\ .
\end{equation}
We call it as the {\bf Hilbert--Stieltjes integral} of the function
$\Phi$ at the point $\tau$.

\bigskip

\begin{proof}
The conclusion of Theorem 2 follows immediately from Lemma 2 and
Corollary 5.
\end{proof} $\Box$

\bigskip

{\bf Corollary 7.} {\it The Hilbert--Stieltjes integral converges
a.e. for every $2\pi -$pe\-rio\-dic bounded function $\Phi:\mathbb
R\to\mathbb R$ that is differentiable a.e.}

\bigskip

{\bf Remark 3.} In particular, the conclusions of Theorem 2 as well
as Corollary 7 hold for every $2\pi -$periodic bounded function of
the class $\mathcal{CBV}(\mathbb R)$.

\bigskip

Of course, Lemmas 1--2, Theorems 1--2, Corollaries 1--7 and the
definition of the Hilbert--Stieltjes integral are extended in the
natural way to complex valued functions $\Phi$.

\bigskip

\cc
\section{On Schwartz--Stieltjes and Cauchy--Stieltjes integrals}

Given a  $2\pi -$periodic bounded function  $\Phi:\mathbb
R\to\mathbb R$ whose set of points of discontinuity has measure
zero, the {\bf Schwartz--Stieltjes integral}
\begin{equation}\label{eqAS} \mathbb S_{\Phi}(z)\ :=\
\frac{1}{2\pi}\ \int\limits_{-\pi}\limits^{\pi}   \frac{e^{it}
+z}{e^{it} -z}\ d\,\Phi(t)\ ,\ \ \ z\in\mathbb D\ ,
\end{equation}
is well--defined  by the previous sections and the function $\mathbb
S_{\Phi}(z)$ is analytic by the definition of the Riemann--Stieltjes
integral and the Weierstrass theorem, see e.g. Theorem I.1.1 in
\cite{Go}. By Theorem 2 and Corollary 2 we have also the following.

\bigskip

{\bf Corollary 8.} {\it Let $\Phi:\mathbb R\to\mathbb R$ be a  $2\pi
-$periodic bounded function that is differentiable a.e. Then
$\mathbb S_{\Phi}(z)$ has finite angular limit $\Phi^{\prime}(\arg\,
\zeta)+i\cdot\mathbb H_{\Phi}(\arg\zeta)$ as $z\to\zeta$ for a.e.
$\zeta\in\partial\mathbb D$.}

\bigskip

It is clear that the definition (\ref{eqAS}), as well as Corollary
8, is extended in the natural way to the case of the complex valued
functions $\Phi$.

\bigskip

Given a $2\pi -$periodic bounded function  $\Phi:\mathbb R\to\mathbb
C$ whose set of points of discontinuity has measure zero, we see
that the integral
\begin{equation}\label{eqAC} \mathbb C_{\Phi}(z)\ :=\
\frac{1}{2\pi}\ \int\limits_{-\pi}\limits^{\pi} \frac{e^{it}
d\,\Phi(t)}{e^{it} -z}\ ,\ \ \ z\in\mathbb D\ ,
\end{equation}
is also well--defined and we call it by the {\bf Cauchy--Stieltjes
integral}. It is easy to see that $\mathbb C_{\Phi}(z)=\frac{1}{2}\
\mathbb S_{\Phi}(z)$.

\bigskip

{\bf Corollary 9.} {\it Let $\Phi:\mathbb R\to\mathbb C$ be a  $2\pi
-$periodic bounded function that is differentiable a.e. Then
\begin{equation}   \label{eqLPC1}
\lim\limits_{z\to\zeta}\ \mathbb C_{\Phi}(z)\ =\ \frac{1}{2}\ \{\,
\Phi^{\prime}(\arg\, \zeta)\ +\ i\cdot\mathbb H_{\Phi}(\arg\zeta)\}
\end{equation} for a.e. $\zeta\in\partial\mathbb D$ along all nontangential paths in
$\mathbb D$ to the point $\zeta$.}

\bigskip

In this connection, note that the Hilbert-Stieltjes integral can be
described in another way for functions $\Phi$ of bounded variation.

Namely, let us denote by $C(\zeta_0,\varepsilon)$,
$\varepsilon\in(0,1)$, $\zeta_0\in\partial\mathbb D$, the rest of
the unit circle $\partial\mathbb D$ after removing its arc
$A(\zeta_0,\varepsilon):=\{ \zeta\in\partial\mathbb D: |\zeta
-\zeta_0|<\varepsilon\}$ and, setting
$$
{\rm I}_{\Phi}(\zeta_0,\varepsilon)\ =\ \frac{1}{2\pi}\
\int\limits_{C(\zeta_0,\varepsilon)} \frac{\zeta
d\,\Phi_*(\zeta)}{\zeta -\zeta_0}\ ,\ \ \ \zeta_0\in\partial\mathbb
D\ ,\ \ \ \ \ {\mbox{where}}\ \ \ \Phi_*(\zeta):=\Phi(\arg\zeta)\ ,
$$
define the {\bf singular integral of the Cauchy--Stieltjes type}
$$
\mathbb I_{\Phi}(\zeta_0)\ =\ \frac{1}{2\pi}\
\int\limits_{\partial\mathbb D} \frac{\zeta\, d\,\Phi}{\zeta
-\zeta_0}\ ,\ \ \ \zeta_0\in\partial\mathbb D\ ,
$$
as a limit of the integrals ${\rm I}(\zeta_0,\varepsilon)$ as
$\varepsilon\to 0$. By paper \cite{Pr} we have that
\begin{equation}   \label{eqPRIWALOW}
\lim\limits_{z\to\zeta}\ \mathbb C_{\Phi}(z)\ =\ \frac{1}{2}\cdot
\Phi^{\prime}(\arg\, \zeta)\ +\ i\cdot\mathbb I_{\Phi}(\zeta)
\end{equation} for a.e. $\zeta\in\partial\mathbb D$ along all nontangential paths
in $\mathbb D$ to the point $\zeta$. Comparing the relations
(\ref{eqLPC1}) and (\ref{eqPRIWALOW}), we come to the following
conclusion.

\bigskip

{\bf Corollary 10.} {\it Let $\Phi:\mathbb R\to\mathbb C$ be a $2\pi
-$periodic function with bounded variation on $[-\pi,\pi]$. Then for
a.e. $\tau\in[-\pi,\pi]$
\begin{equation}   \label{eqEQ}
\mathbb H_{\Phi}(\tau)\ =\ 2\cdot \mathbb I_{\Phi}(e^{i\tau})\ .
\end{equation}}


\cc
\section{Representation of the Luzin construction}

The following deep (non--trivial) result of Luzin was one of the
main theorems of his dissertation, see e.g. his paper \cite{L1},
dissertation \cite{L2}, p. 35, and its reprint \cite{L}, p. 78,
where one may assume that $\Phi(0)=\Phi(1)=0$, cf. also \cite{RY}.

\medskip

{\bf Theorem A.} {\it\, For any measurable function $\varphi :
[0,1]\to \mathbb R$, there is a continuous function $\Phi : [0,1]\to
\mathbb R$ such that $\Phi^{\prime}=\varphi$ a.e.}

\medskip

Just on the basis of Theorem A, Luzin proved the next significant
result of his dissertation, see e.g. \cite{L}, p. 80, that was key
to establish the corresponding result on the boundary value Hilbert
problem for analytic functions in \cite{R1}.

\medskip

{\bf Theorem B.} {\it\, Let $\varphi(\vartheta)$ be real,
measurable, almost everywhere finite and have the period $2\pi$.
Then there exists a harmonic function $U$ in the unit disk $\mathbb
D$ such that $U(z)\to \varphi(\vartheta)$ for a.e. $\vartheta$ as
$z\to e^{i\vartheta}$ along any nontangential path.}

\medskip

Note that the Luzin dissertation was published in Russian as the
book \cite{L} with comments of his pupils Bari and Men'shov  only
after his death. A part of its results was also printed in  Italian
\cite{Lu}. However, Theorem A was published with a complete proof in
English in the book \cite{S} as  Theorem VII(2.3). Hence Frederick
Gehring in \cite{Ge} has rediscovered Theorem B and his proof on the
basis of Theorem A has in fact coincided with the original proof of
Luzin. Since the proof is very short and nice and has a common
interest, we give it for completeness here.

\medskip

\begin{proof}
By Theorem A  we can find a continuous function $\Phi(\vartheta)$
such that $\Phi^{\prime}(\vartheta)=\varphi(\vartheta)$ for a.e.
$\vartheta$. Considering the Poisson integral
$$
u(re^{i\vartheta})\ =\ \frac{1}{2\pi}\
\int\limits_{0}\limits^{2\pi}\frac{1-r^2}{1-2r\cos(\vartheta-t)+r^2}\
\Phi(t)\ dt
$$
for $0<r<1$, $u(0):=0$, we see by the Fatou result, see e.g. 3.441
in \cite{Z}, p. 53, or Theorem IX.1.2 in \cite{Go}, that
$\frac{\partial}{\partial\vartheta}\ u(z)\to
\Phi^{\prime}(\vartheta)$ as $z\to e^{i\vartheta}$ along any
nontangential path whenever $\Phi^{\prime}(\vartheta)$ exists. Thus,
the conclusion follows for the function $U(z)\ =\
\frac{\partial}{\partial\vartheta}\ u(z)$. \end{proof} $\Box$

\medskip

{\bf Remark 4.} Note that the given function $U$ is harmonic in the
punctured unit disk $\mathbb D\setminus\{0\}$ because the function
$u$ is harmonic in $\mathbb D$ and the differential operator
$\frac{\partial}{\partial\vartheta}$ is commutative with the Laplace
operator $\Delta$. Setting $U(0)=0$, we see that
$$
U(re^{i\vartheta})\ =\ -\frac{r}{\pi}\
\int\limits_{0}\limits^{2\pi}\frac{(1-r^2)\sin(\vartheta-t)}{(1-2r\cos(\vartheta-t)+r^2)^2}\
\Phi(t)\ dt\ \to\ 0\ \ \ \ \ \ \ \mbox{as}\ r\to 0\ ,
$$
i.e. $U(z)\to U(0)$ as $z\to 0$, and, moreover, the integral of $U$
over each circle $|z|=r$, $0<r<1$, is equal to zero. Thus, by the
criterion for a harmonic function on the averages over circles we
have that $U$  is harmonic in $\mathbb D$. The alternative argument
for the latter is the removability of isolated singularities for
harmonic functions, see e.g. \cite{N}.

\medskip

In the connection with Theorem B, it is necessary also to mention
the paper \cite{M} contained the Men'shoff theorem on the existence
of a trigonometric series that is  a.e. convergent to a prescribed
measurable function $\varphi : (0,2\pi)\to\mathbb R$.

\medskip

Corollary 5.1 in \cite{R1} has strengthened Theorem B as the next,
see also \cite{R2}.

\medskip

{\bf Theorem C.} {\it\, For each (Lebesgue) measurable function
$\varphi :\partial\mathbb D\to \mathbb R$, the space of all harmonic
functions $u:\mathbb D\to\mathbb R$ with the angular limits
$\varphi(\zeta)$ for a.e. $\zeta\in\partial\mathbb D$ has the
infinite dimension.}

\medskip

{\bf Remark 5.} One can find in \cite{RY} more refined results which
are counterparts of Theorems A, B and C in terms of logarithmic
capacity that makes possible to extend the theory of boundary value
problems to the so--called $A-$harmonic functions corresponding to
generalizations of the Laplace equation in inhomogeneous and
anisotropic media, see also \cite{Y}.

By the well--known Lindel\"of maximum principle, see e.g. Lemma 1.1
in \cite{GM}, it follows the uniqueness theorem for the Dirichlet
problem in the class of bounded harmonic functions $u$ on the unit
disk $\mathbb D = \{ z\in\mathbb{C}: |z|<1\}$. In general there is
no uniqueness theorem in the Dirichlet problem for the Laplace
equation even under zero boundary data.

\medskip

Many such nontrivial solutions for the Laplace equation can be given
just by the Poisson-Stieltjes integral
\begin{equation}\label{eqPS}
\mathbb U_{\Phi}(z)\ :=\ \frac{1}{2\pi}\
\int\limits_{0}\limits^{2\pi}P_r(\vartheta -t)\ d\Phi(t)\ ,\ \ \
z=re^{i\vartheta}, \ r<1\ ,
 \end{equation}
with an arbitrary {\bf singular function} $\Phi:[0,2\pi]\to\mathbb
R$, i.e., where $\Phi$ is of bounded variation and $\Phi^{\prime}=0$
a.e. Indeed, by the Fatou theorem, see e.g. Theorem I.D.3.1 in
\cite{Ko}, $\mathbb U_{\Phi}(z)\to\Phi^{\prime}(\Theta)$ as $z\to
e^{i\Theta}$ along any nontangential path whenever
$\Phi^{\prime}(\Theta)$ exists, i.e. $\mathbb U_{\Phi}(z)\to0$ as
$z\to e^{i\Theta}$ for a.e. $\Theta\in[0,2\pi]$ along any
nontangential paths for every singular function $\Phi$.

\medskip

{\bf Example 2.} {\it The first natural example is given by the
formula (\ref{eqPS}) with $\Phi(t)=\varphi(t/2\pi)$ where
$\varphi:[0,1]\to[0,1]$ is the well--known {\bf Cantor function},}
see e.g. \cite{DMRV} and further references therein.

\medskip

{\bf Example 3.} {\it However, the simplest example of such a kind
is just
$$
u(z)\ :=\ P_r(\vartheta -\vartheta_0)\ =
\frac{1-r^2}{1-2r\cos(\vartheta -\vartheta_0) +r^2}\ ,\ \ \
z=re^{i\vartheta}, \ r<1\ .
$$
We see that $u(z)\to0$ as $z\to e^{i\Theta}$ for all
$\Theta\in(0,2\pi)$ except $\Theta=\vartheta_0$.}

\bigskip

The construction of Luzin can be described as the Poisson--Stieltjes
integral.

\bigskip

{\bf Theorem 3.} {\it The harmonic function $U$ in Theorem B has the
representation
\begin{equation}\label{eqPS1}
\mathbb U_{\Phi}(z)\ =\frac{1}{2\pi}\
\int\limits_{-\pi}\limits^{\pi} P_r(\vartheta -t)\ d\Phi(t)\ \ \ \
\forall\ z=re^{i\vartheta}, \ r\in(0,1)\ ,\ \vartheta\in[-\pi,\pi]\
,
 \end{equation}
where $\Phi:[-\pi,\pi]\to\Bbb R$ is the continuous Luzin function
with $\Phi^{\prime}=\varphi$ a.e.}

\bigskip

{\bf Corollary 11.} The conjugate harmonic function $\mathbb
V_{\Phi}$ has finite angular limits
$$
\lim\limits_{z\to\zeta}\mathbb V_{\Phi}(z)\ =\ \mathbb H_{\Phi}(\arg
\zeta)\ \ \ \ \ \ \ \ {\mbox{for\ a.e.}}\ \ \
\zeta\in\partial\mathbb D\ .
$$

\begin{proof} Indeed, choosing in Proposition 1 $g(t)=\Phi(t)$ and $f(t)=P_r(\vartheta
-t)$, $t\in[-\pi,\pi]$, for every fixed $z=re^{i\vartheta}$, $r<1$,
$\vartheta\in[-\pi,\pi]$, we obtain that
\begin{equation}\label{eqPS2INVERSE}
\int\limits_{-\pi}\limits^{\pi} \Phi(t)\cdot
\frac{\partial}{\partial\vartheta}\, P_r(\vartheta -t)\ dt\ =\
\int\limits_{-\pi}\limits^{\pi} P_r(\vartheta -t)\ d\Phi(t)\ \ \ \
\forall\  r\in(0,1)\ ,\ \vartheta\in[-\pi,\pi]
 \end{equation}
because by the $2\pi -$periodicity of the given function $g$ and $f$
the right hand side in  (\ref{eqPARTS}) is equal to zero, $f\in
C^1$, and
$$
\frac{\partial}{\partial\vartheta}\, P_r(\vartheta -t)\ =\ -\
\frac{\partial}{\partial t}\, P_r(\vartheta -t)\ \ \ \ \forall\
r\in(0,1)\ ,\ \vartheta\in[-\pi,\pi]\ ,\ t\in[-\pi,\pi]\ .
$$
The relation (\ref{eqPS1}) follows from (\ref{eqPS2INVERSE}) because
it was in the proof of Theorem B
$$
U(z)\ =\ \frac{1}{2\pi}\  \int\limits_{-\pi}\limits^{\pi}
\Phi(t)\cdot \frac{\partial}{\partial\vartheta}\, P_r(\vartheta -t)\
dt\ \ \ \ \ \forall\ z=re^{i\vartheta}, \  r\in(0,1)\ ,\
\vartheta\in[-\pi,\pi]\ .
$$\end{proof} $\Box$

\bigskip \noindent
{\bf Vladimir Il'ich Ryazanov,}\\
Institute of Applied Mathematics and Mechanics,\\
National Academy of Sciences of Ukraine, room 417,\\
19 General Batyuk Str., Slavyansk, 84116, Ukraine,\\
vlryazanov1@rambler.ru , vl.ryazanov1@gmail.com

\end{document}